\documentclass[11pt]{amsart}

\usepackage{amsthm}
\usepackage{amsmath}
\usepackage{amssymb}
\usepackage{color}

\newtheorem{thm}{Theorem}[section]
\newtheorem{theorem}[thm]{Theorem}

\newtheorem{corollary}[thm]{Corollary}
\newtheorem{problem}[thm]{Problem}

\newtheorem{example}[thm]{Example}

\newtheorem{definition}[thm]{Definition}
\theoremstyle{remark}

\newcommand{\RR}{\mathbb R}

\newcommand{\CC}{\mathbb C}
\newcommand{\HH}{\mathbb H}

\begin{document}

\title{an elementary proof of a fundamental result in real phase retrieval}
\author[Casazza, Tremain
 ]{Peter G. Casazza and Janet C. Tremain}
\address{Department of Mathematics, University
of Missouri, Columbia, MO 65211-4100}

\thanks{The first author was supported by
 NSF DMS 1906725}

\email{Casazzap@missouri.edu; tremainjc@missouri.edu}


\begin{abstract}
Edidin \cite{E} proved a fundamental result in phase retrieval:  Theorem:  A family of orthogonal projections $\{P_i\}_{i=1}^m$
does phase retrieval in $\RR^n$ if and only if for every $0\not= x\in \RR^n$, the family $\{P_ix\}_{i=1}^m$ spans $\RR^n$.
The proof of this result relies on Algebraic Geometry and so is inaccessible to many people in the field. We will give an
elementary proof of this result using only frame theory.  We will also solve the complex version of this result by showing
that the "if" part fails and the "only if" part holds in $\CC^n$. Finally, we will show that these techniques can be used
to verify two classifications of norm retrieval.
\end{abstract}

\maketitle

\section{Introduction}
We start with frames.

\begin{definition}
A family of vectors $\{x_i\}_{i=1}^m$ is a {\bf frame} for $\HH^n$ if there are constants $0<A\le B<\infty$ so that
for every $x\in \HH^n$:
\[ A \|x\|^2 \le \sum_{i=1}^m|\langle x,x_i\rangle|^2 \le B\|x\|^2.\]
\end{definition}

We have
\begin{enumerate}
\item $A,B$ are a set of {\bf lower} (respectively, {\bf upper}) frame bounds.
\item If $A=B$ this is an {\bf A-tight frame} and if $A=B=1$ it is a {\bf Parseval frame}.
\item The operator $T:\HH^n\rightarrow \ell_2(m)$ given by 
\[Tx=(\langle x,x_1\rangle, \langle x,x_2\rangle, \cdots, \langle x,x_m\rangle)\]
is called the {\bf analysis operator} of th{\bf frame operator}. The frame operator is a apositive, invertible, self-adjoint operator on the space.
\end{enumerate}

Phase retrieval for Hilbert space frames was introduced in \cite{BCE} and quickly became one of the most active areas of
research in frame theory. 

\begin{definition}
A frame $\{x_i\}_{i=1}^m$ does {\bf phase retrieval} in $H^n$ if whenever $x,y\in \HH^n$ and
\[ |\langle x,x_i\rangle|=|\langle y,x_i\rangle|,\mbox{ for all }i=1,2,\ldots,m,\]
then $x=cy$ for some $|c|=1$.
\end{definition}

 In some applications a signal must be reconstructed from the norms of higher dimensional components. In X-ray crystallography
such a problem arises in crystal twinning \cite{Cr}. In this scenario, there exists a similar phase retrieval problem \cite{CCJ}.

\begin{definition}
A family of subspaces $\{W_i\}_{i=1}^m$ with respective projections $\{P_i\}_{i=1}^m$ does {\bf phase retrieval} if whenever
$x,y\in \HH^n$ and
\[ \|P_ix\|=\|P_iy\|,\mbox{ for all }i=1,2,\ldots,m,\]
then $x=cy$ for some $|c|=1$.
\end{definition}

\section{Background}

We need:

\begin{definition}
A family of vectors $\{x_i\}_{i=1}^n$ in $\RR^n$ or $C^n$ has the {\bf complement property} if whenever we partition the set of vectors into
two subsets, one of them must span the space.
\end{definition}

A fundemantal result here appeared in \cite{BCE}:

\begin{thm}
A family of vectors $\{x_i\}_{i=1}^m$ in $\RR^n$ does phase retrieval if and only if it has the complement property.
\end{thm}

\begin{corollary}
If a family of vectors $\{x_i\}_{i=1}^m$ satisfies phase retrieval in $\RR^n$, then $m\ge 2n-1$.
\end{corollary}

It is easy to get families of vectors doing phase retrieval.

\begin{definition}
A frame $\{x_i\}_{i=1}^m$ is {\bf full spark} in $\HH^n$ if for every $I\subset [m]$ with $|I|=n$, span $\{x_i\}_{i\in I}=\HH^n$.
\end{definition}

It follows that any full spark set of $\ge 2n-1$ vectors does phase retrieval. There are full spark sets with an arbitrary number of
vectors in $\HH^n$. Let $\{e_i\}_{i=1}^n$ be an orthonormal basis for $\HH^n$. For every $I\subset [n]$ with $|I|=n-1$, 
span $\{x_i\}_{i\in I}$ is a hyperplane. Since we only have a finite number of hyperplanes, there is a vector $x_{m+1}\in \HH^n$
so that $x_{m+1}\notin span\ \{x_i\}_{i\in I}$ for any $I\subset [n]$ with $|I|=n-1$, and  hence
 $\{x_i\}_{i=1}^{m+1}$ is full spark. Now continue this. Choose a vector $x_{m+2}\in \HH^n$ so that $x_{m+2}\in \HH^n$
and $x_{m+2}\notin span\ \{x_i\}_{i\in I}$ for all $I\subset [m+1]$ and $|I|=n-1$. Hence, $\{x_i\}_{i=1}^{n+2}$ is full spark.

The complex case is much more complicated. It is known \cite{BCMN} that phase retrieval can be done in $\CC^n$ with
$4n-4$ vectors. It was believed this result was best possible. But Vinzant \cite{V} showed that phase retrieval can be
done in $\CC^4$ with 11 vectors (in other words $4n-5$ vectors). So the major open problem here is:

\begin{problem}
How many vectors in $C^n$ does it take to do phase retrieval?
\end{problem}

For phase retrieval by projections the problem is still more complicated.
In \cite{CCJ} it was shown:

\begin{thm}
Let $\{P_i\}_{i=1}^n$ be orthogonal projections in $\RR^n$ with $P_i(\RR^n)=W_i$, for $i=1,2,\ldots,m$. The following
are equivalent:
\begin{enumerate}
\item $\{W_i\}_{i=1}^m$ allows phase retrieval.
\item Whenever $\{x_{ij}\}_{j=1}^{I_i}$ is an orthonormal basis for $W_i$, the family $\{x_{ij}\}_{i=1,j=1}^{ m\ \ \ I_i}$
does phase retrieval.
\end{enumerate}
\end{thm}

Recently, \cite{CG} the following was shown:

\begin{theorem}
It takes at least $2n-2$ subspaces to do phase retrieval in $\RR^n$.
\end{theorem}

The following theorem in \cite{X} shows a particular case of $\RR^4$ in which phase retrieval is possible with $2n-2= 2(4)-2=6$ subspaces.

\vspace{0.1 in}
\begin{theorem}
There are 6 two-dimensional subspaces of $\RR^4$ satisfying phase retrieval.
\end{theorem}

In \cite{A} it was shown that there are 6 hyperplanes in $\RR^4$ doing phase retrieval.

Another open problem in this area is motivated by this theorem of 
Edidin \cite{E}:
\vspace{0.1 in}

\begin{theorem}
If $n=2^k+1$, for some $k$, then it takes $2n-1$ subspaces of $\RR^n$ to do phase retrieval.
\end{theorem}

So now the main open problem in this area is:
\vspace{0.1 in}

\begin{problem}
Can we do phase retrieval in $\RR^n$ with $2n-2$ subspaces whenever $n\not= 2^k+1$ for any $k$.
\end{problem}

Again, the complex case is  much more complicated.

\begin{problem}
What is the minimum number of projections required to do phase retrieval in $\CC^n$?
\end{problem}


\vspace{0.1 in}

Edidin \cite{E} proved a fundamental result in phase retrieval: 

\begin{thm} \label{T}A family of orthogonal projections $\{P_i\}_{i=1}^m$
does phase retrieval in $\RR^n$ if and only if for every $0\not= x\in \RR^n$, the family $\{P_ix\}_{i=1}^m$ spans $\RR^n$.
\end{thm}
Because the proof relies on Algebraic Geometry, it has been inaccessible to some people in the field. We will give an
elementary proof of this result without Algebraic Geometry. We will also resolve the complex case by showing that
the "if" part of the theorem fails and the "only if" part of the theorem holds in $C^n$. We will also show that these methods
verify several classifications of norm retrieval.
\vskip10pt

\noindent {\bf Remark}: Edidin's Theorem also holds for symmetric matricies $\{A_i\}_{i=1}^m$ instead of projections
\cite{X}

\section{Phase Retrieval}
We start with the real case.

\subsection{Proof of Theorem \ref{T} for the real case.}
\begin{proof}
$\Rightarrow$: We do the contrapositive.
Assume there is an $x\in R^n$with $\|x\|=1$ so that $span \{P_ix\}_{i=1}^m \not=  R^n$. Choose $y\in R^n$ with
$\|y\|=1$ and $y\perp span\ \{P_ix\}_{i=1}^m$. Now,
\[ 0=\langle y,P_ix\rangle = \langle P_iy,P_ix\rangle,\mbox{ for all }i=1,2,\ldots,m.\]
Since $P_ix\perp P_iy$,
\[ \|P_i(x+y)\|^2=\|P_ix\|^2+\|P_iy\|^2 = \|P_i(x-y)\|^2.\]
But, $x+y\not= \pm(x-y)$. So $\{P_i\}_{i=1}^m$ fails phase retrieval.

\vskip12pt
$\Leftarrow $:
Let $W_i = P_i(R^n)$. If $\{W_i\}_{i=1}^m$ fails phase retrieval, then there exists orthonormal bases $\{x_{ij}\}_{j=1}^{I_i}$
for $W_i$ which can be divided into two non-spanning sets $\{x_{ij}\}_{(i,j)\in I}$ and $\{x_{ij}\}_{(i,j)\in I^c}$. Choose
\[ 0\not= x \perp span\ \{x_{ij}\}_{(i,j)\in I}.\]
Then, $P_ix \in span\ \{x_{ij}\}_{(i,j)\in (I^c\cap W_i)}$ for all i. But this family of sets does not span.

\end{proof}

\subsection{Proof of the Complex Case}

\begin{proof}
$\Rightarrow$:  Assume there is an $x\in \CC^n$ with $\|x\|=1$ so that $\{P_ix\}_{i=1}^m$ does not span.  Choose
$y\in \CC^n$ with $\|y\|=1$ so that $y\perp P_ix$, for all $i=1,2,\ldots,m$.  Now
\begin{align*}
\|P_i(x+y)\|^2 &= \|P_ix\|^2+\|P_iy\|^2+2Re\langle P_ix,P_iy\rangle \\
&= \|P_ix\|^2+\|P_iy\|^2+2 Re\langle P_ix,y\rangle\\
&= \|P_ix\|^2+\|P_iy\|^2 .
\end{align*}
Similarly,
\[ \|P_i(x-y)\|^2 = \|P_ix\|^2+\|P_iy\|^2.\]
So if $\{P_i\}_{i=1}^m$ does phase retrieval, then there is some $|c|=1$ so that $x+y=c(x-y)$. If $c=1$ then
$x+y=x-y$ and so $y=0$, a contradiction.  Otherwise, $(1+c)y=(c-1)x$ and so
$x=\frac{1+c}{c-1}y$. Thus, $P_ix=\frac{1+c}{c-1}P_iy$. But this is impossible since
$x,y\not= 0$ and
\[ \langle P_ix,P_iy\rangle = \langle P_ix,y\rangle =0, \mbox{ for all }i=1,2,\ldots,m.\]
\end{proof}

It is known (see for example, \cite{V}) that for a family of vectors $\{x_i\}_{i=1}^m$ to do phase retrieval in $C^n$, m must be on the order of
4n.

\begin{example}
Let $\{x_i\}_{i=1}^{2n-1}$ be a full spark family of vectors in $\CC^n$. This family cannot do phase retrival. But if $0\not= x\in \CC^n$,
then $|\{1\le i \le 2n-1:\langle x,x_i\rangle \not= 0\}|\ge n$.  For otherwise, $|\{1\le i \le 2n-1:\langle x,x_i\rangle =0\}|\ge n$
and since this family is full spark, $x=0$ a contradiction. By the first inequality and the fact that this family is full spark,
$\{x_i:\langle x,x_i\rangle \not= 0\}$ spans the space.
\end{example}

\noindent {\bf Remark}: There is a complete complex version of Edidin's theorem: See Theorem 2.1 (3) in \cite{X}.
\vskip10pt

A completely different approach to a complex version of Edidin's theorem was proved in \cite{CC}. In this paper they associate vectors
in $\CC^n$  with
rank 2 projections in $\RR^{2n}$.

\section{Norm Retrieval}

Norm retrieval is a weakening of phase retrieval.

\begin{definition}
A frame $\{x_i\}_{i=1}^m$ in $\HH^n$ does {\bf norm retrieval} if whenever $x,y\in \HH^n$ and
\[ |\langle x,x_i\rangle|=|\langle y,x_i\rangle|,\mbox{ for all }i=1,2,\ldots,m,\]
then $\|x\|=\|y\|$. 
\end{definition}

If $\{x_i\}_{i=1}^m$ does phase retrieval, then it does norm retrieval. But an orthonormal basis $\{e_i\}_{i=1}^n$ does norm
retrieval in $\HH^n$ and fails phase retrieval.

The techniques from the previous section classify norm retrieval \cite{F}.

\begin{theorem}
A family of vectors $\{x_i\}_{i=1}^m$ does norm retrieval if and only if whenever $I\subset [m]$ and $x\perp x_i$,
for all $i\in I$ and $y\perp x_i$ for all $i\in I^c$, then $x\perp y$.
\end{theorem}

\begin{proof}
$\Leftarrow$: Assume $|\langle x,x_i\rangle| = |\langle y,x_i\rangle|$, for all $i=1,2,\ldots,m$. Let
\[ I=\{1\le i \le m:\langle x,x_i\rangle = \langle y,x_i\rangle\}.\]
Then, $x+y\perp x_i$ for all $i\in I^c$ and $x-y\perp x_i$ for all $\in I$. By assumption, $x+y\perp x-y$. Now,
\[ 0=\langle x+y,x-y\rangle = \|x\|^2-\|y\|^2,\mbox{ and so } \|x\|=\|y\|.\]
So $\{x_i\}_{i=1}^m$ does norm retrieval.
\vskip10pt

$\Rightarrow$: Assume $\{x_i\}_{i=1}^m$ does norm
retrieval and assume $I\subset [m]$ and $x\perp x_i$ for all $i\in I$ and $y\perp x_i$ for all $i\in I^c$.  Then
\[ |\langle x+y,x_i\rangle|=|\langle x-y,x_i\rangle|, \mbox{ for all } i=1,2,\ldots,m.\]
Since $\{x_i\}_{i=1}^m$ does norm retrieval, 
\begin{align*}
\|x+y\|^2 &= \|x\|^2+\|y\|^2+2\langle x,y\rangle\\
&= \|x-y\|^2\\
&= \|x\|^2+\|y\|^2-2\langle x,y\rangle.
\end{align*}
Hence, $\langle x,y\rangle =0$ and $x\perp y$.
\end{proof}

There is a corresponding notion for projections.

\begin{definition}
A family of projections does {\bf norm retrieval} if whenever $x,y\in \HH^n$ and
\[ \|P_ix\|=\|P_iy\|, \mbox{ for all }i=1,2,\ldots,m,\]
then $\|x\|=\|y\|$.
\end{definition}

In general, if $\{P_i\}_{i=1}^m$ does phase retrieval, it may not be true that $\{(I-P_i)\}_{i=1}^m$ does phase retrieval.
For example, if \cite{A} $\{W_i\}_{i=1}^6$ are hyperplanes doing phase retrieval in $\RR^4$, then $\{W_i^{\perp}\}_{i=1}^6$ are
only 6 vectors in $\RR^4$ and hence cannot do phase retrieval. What is needed here is norm retrieval.

\begin{theorem}
If $\{P_i\}_{i=1}^m$ does phase retrieval, then $(I-P_i)\}_{i=1}^m$ does phase retrieval in $\RR^n$ if and only if $\{(I-P_i)\}_{i=1}^m$ does
norm retrieval.
\end{theorem}

\begin{proof} $\Rightarrow$:
If $\{(I-P_i)\}_{i=1}^m$ does phase retrieval, it does norm retrieval.
\vskip10pt
$\Leftarrow$:
Assume $\{P_i\}_{i=1}^m$ does phase retrieval and $\{(I-P_i)\}_{i=1}^m$ does norm retrieval. If $x,y\in \HH^n $ and
\[ \|(I-P_i)x\|=\|(I-P_i)y\|,\mbox{ for all }i=1,2,\ldots,m,\]
then by assumption, $\|x\|=\|y\|$. Also,
\[ \|(I-P_i)x\|^2=\|x\|^2-\|P_ix\|^2=\|(I-P_i)y\|^2=\|y\|^2-\|P_iy\|^2.\]
Hence, $\|P_ix\|=\|P_iy\|$, for all $i=1,2,\ldots,m$ and since $\{P_i\}_{i=1}^m$ does phase retrieval, $x=\pm y$. I.e. $\{(I-P_i)\}_{i=1}^m$ does phase retrieval.
\end{proof}

There is a form of norm retrieval for projections.

\begin{theorem}
Let $\{W_i\}_{i=1}^m$ be subspaces of $\HH^n$ with corresponding projections $\{P_i\}_{i=1}^m$. The following are
equivalent:
\begin{enumerate}
\item $\{W_i\}_{i=1}^m$ does norm retrieval.
\item For every orthonormal basis $\{x_{ij}\}_{j=1}^{I_i}$ of $W_i$ the family $\{x_{ij}\}_{i=1,j=1}^{\ m\ I_i}$ does norm
retrieval.
\end{enumerate}
\end{theorem}

\begin{proof}
$(1) \Rightarrow (2)$: Assume $\{W_i\}_{i=1}^m$ does norm retrieval. Let $\{x_{ij}\}_{j=1}^{I_i}$ be an orthonormal basis of
$W_i$. Assume
\[ |\langle x,x_{ij}\rangle|=|\langle y,x_{ij}\rangle|,\mbox{ for all }(i,j).\]
Then, for all $i=1,2,\ldots,m$,
\[ \|P_ix\|^2=\sum_{j=1}^{I_i}|\langle x,x_{ij}\rangle |=\sum_{j=1}^{I_i}|\langle y,x_{ij}\rangle = \|P_iy\|^2.\]
Since $\{P_i\}_{i=1}^m$ does norm retrieval, $\|x\|=\|y\|$. So $\{x_{ij}\}_{i=1,j=1}^{  m\ \ \ I_i}$ does norm retrieval.
\vskip10pt
$(2) \Rightarrow (1)$: Assume (2). Let $x,y\in \HH^n$ satisfy $\|P_ix\|=\|P_iy\|$, for all $i=1,2,\ldots,m$.
If $P_ix=\pm P_iy$ then $|\langle x,x_{ij}\rangle|=|\langle y,x_{ij}\rangle|$, for all $j=1,2,\ldots,I_i$.
Otherwise, let
\[ x_{i1}:= \frac{P_ix+P_iy}{\|P_ix+P_iy\|},\ \ x_{i2}= \frac{P_ix-P_iy}{\|P_ix-P_iy\|}.\]
Then
\[ \langle x_{i1},x_{i,2}\rangle = \frac{1}{\|P_ix+P_iy\|\|P_ix-P_iy\|}(\|P_ix\|^2-\|P_iy\|^2 +\langle P_ix,P_iy\rangle -\langle
P_iy,P_ix\rangle)=0.\]
So we can extend $\{x_{i1},x_{i2}\}$ to an orthonormal basis for $W_i$, $\{x_{ij}\}_{j=1}^{I_i}$. Note that $P_ix,P_iy \in 
span\{x_{i1},x_{i2}\}$. So,
\[ \langle x,x_{ij}\rangle = \langle y,x_{ij}\rangle=0,\mbox{ for all } j=3,4,\ldots,I_i.\]
We also have by the above
\begin{eqnarray*}
|\langle x,P_ix+P_iy\rangle|&=|\|P_ix\|^2+\langle P_ix,P_iy\rangle|\\
&= |\|P_iy\|^2+\langle P_ix,P_iy\rangle|\\
&= |\langle P_iy,P_iy+P_ix\rangle|\\
&= |\langle y,P_ix+P_iy\rangle|.
\end{eqnarray*} 
Similarly, $|\langle x,P_ix-P_iy\rangle|=|\langle y,P_ix-P_iy\rangle|$. Hence,
\[ |\langle x,x_{ij}\rangle|=|\langle y,x_{ij}\rangle|, \mbox{ for all }j=1,2,\ldots,I_i.\]
By (2), $\|x\|=\|y\|$, and so $\{P_i\}_{i=1}^m$ does norm retrieval.
\end{proof}

We have an Edidin type theorem for norm retrieval \cite{CC}.

\begin{theorem}
Given projections $\{P_i\}_{i=1}^m$ in $\RR^n$, the following are equivalent:
\begin{enumerate}
\item $\{P_i\}_{i=1}^m$ does norm retrieval.
\item For every $0\not= x \in \RR^n$, we have
\[ \left [ span\ \{P_ix\}\right ] ^{\perp} \subset x^{\perp}.\]
\item For every $x\in \HH^n$ we have
$x\in span\ \{P_ix\}_{i=1}^m$.
\end{enumerate}
 
\end{theorem}

\begin{proof}
$(2)\Rightarrow (1)$: We will prove the contrapositive. If norm retrieval fails, then there are vectors $x,y\in \RR^n$ with
$\|P_ix\|=\|P_iy\|$, for all $i=1,2,\ldots,m$ but $\|x\|\not= \|y\|$. This implies if $v=x+y$ and $w=x-y$, then $v,w$ are not
orthogonal. But, $w\in [span\ \{P_iv\}_{i=1}^m]^{\perp}$. So (2) fails.
\vskip10pt
$(2)\Leftrightarrow (3)$: This is immediate.
\vskip10pt
$(2)\Rightarrow (1)$: Again by the contrapositive. Assume there exists $v,w\in \RR^n$ which are not orthogonal (so $w\not= 0$) but
$w\perp span\ \{P_iv\}_{i=1}^m$ and so $P_iw \perp P_iv$, for all $i=1,2,\ldots,m.$ Let $v=x+y$ and $w=x-y$. Since $v,w$ are not
orthogonal, $\|x\|\not= \|y\|$. But,
\begin{align*}
\|P_i(v+w)\|^2&= \langle P_iv+P_iw,P_iv+P_iw\rangle\\
&= \|P_iv\|^2+\|P_iw\|^2\\
&= \langle P_iv-P_iw,P_iv-P_iw\rangle\\
&= \|P_i(v-w)\|^2.
\end{align*} 
So $\{P_i\}_{i=1}^m$ fails norm retrieval.
\end{proof}

The projection form of this is:

\begin{theorem}
Let $\{W_i\}_{i=1}^m$ be subspaces of $\HH^n$ with respective projections $\{P_i\}_{i=1}^m$. The following are equivalent:
\begin{enumerate}
\item $\{W_i\}_{i=1}^m$ does norm retrival.
\item If $\{x_{ij}\}_{j=1}^{I_i}$ is an orthonormal basis for $W_i$, and $S\subset \{(ij):1\le i \le m,\ 1\le i \le I_i\}$ and
$x,y\in \HH^n$ satisfy 
\[ x\perp x_{ij}\mbox{ for }(ij)\in S \mbox{ and } y\perp x_{ij}\mbox{ for } (ij)\in S^c,\]
then $x \perp y$.
\end{enumerate}
\end{theorem}

\begin{proof}
$(1) \Rightarrow (2)$: Assume $\{W_i\}_{i=1}^m$ does norm
retrieval and assume $S\subset \{(i,j):1\le i \le m,\ 1\le j \le I_i\}$,  and $x\perp x_{ij}$ for all $(i,j)\in S$ 
and $y\perp x_{iij}$ for all $i(i,j)\in S^c$.  Then
\[ |\langle x+y,x_{ij}\rangle|=|\langle x-y,x_{ij}\rangle|, \mbox{ for all } (i,j).\]
Since $\{W_i\}_{i=1}^m$ does norm retrieval, 
\begin{align*}
\|x+y\|^2 &= \|x\|^2+\|y\|^2+2\langle x,y\rangle\\
&= \|x-y\|^2\\
&= \|x\|^2+\|y\|^2-2\langle x,y\rangle.
\end{align*}
Hence, $\langle x,y\rangle =0$ and $x\perp y$.
\vskip10pt
$(2) \Rightarrow (1)$: Assume $|\langle x,x_{ij}\rangle|=|\langle y,x_{ij}\rangle|$ for all $(i,j)$. Let
\[ S=\{(i,j):\langle x,x_{ij}\rangle = \langle y,x_{ij}\rangle\}.\]
Then $x+y \perp x_{ij}$ for all $(i,j)\in S^c$ and $x-y \perp x_{ij}$ for all $(i,j)\in S$. By assumpgion, 
$x+y\perp x-y$ so
\[ 0=\langle x+y,x-y\rangle = \|x\|^2-\|y\|^2,\mbox{ and so } \|x\|=\|y\|.\]
So $\{x_{ij}\}$ does norm retrieval.
\end{proof}

\end{document}